\documentclass{IEEEtran4PSCC}

\usepackage{cite}

\ifCLASSINFOpdf
   \usepackage[pdftex]{graphicx}
\else
   \usepackage[dvips]{graphicx}
\fi
\usepackage[cmex10]{amsmath}
\usepackage{algorithmic}

\usepackage{array}

\ifCLASSOPTIONcompsoc
 \usepackage[caption=false,font=normalsize,labelfont=sf,textfont=sf]{subfig}
\else
 \usepackage[caption=false,font=footnotesize]{subfig}
\fi
\usepackage{stfloats}
\usepackage{url}

\usepackage{mathrsfs}
\usepackage{amssymb,amsfonts}
\usepackage{color}
\usepackage[normalem]{ulem}
\usepackage{enumitem}
\usepackage[utf8]{inputenc}
\usepackage{textcomp}
\usepackage{eurosym}
\usepackage{accents}
\newcommand{\ubar}[1]{\underaccent{\bar}{#1}}
\usepackage{hyperref}
\usepackage{blindtext}
\usepackage{xcolor}
\usepackage[normalem]{ulem}
\usepackage{booktabs}
\usepackage{acronym}
\usepackage[mathscr]{eucal}
\usepackage{mathtools}
\usepackage{graphicx}
\usepackage{pifont}
\newcommand{\xmark}{\ding{56}}%
\usepackage{multirow}
\usepackage{xpatch}
\usepackage{nomencl}

\newcommand{\nomunit}[1]{%
\renewcommand{\nomentryend}{\hspace*{\fill}#1}}
\makenomenclature

\makeatletter
\patchcmd{\thenomenclature}{\section*}{\section}{}{}
\makeatother

\usepackage{etoolbox}
\renewcommand\nomgroup[1]{%
  \item[\bfseries
  \ifstrequal{#1}{I}{Indexes and Sets}{%
  \ifstrequal{#1}{P}{Parameters}{%
  \ifstrequal{#1}{V}{Variables}{}}}%
]}

\acrodef{dam}[DAM]{Day-Ahead Market}
\acrodef{idm}[IDM]{Intra-Day Market}
\acrodef{rtm}[RTM]{Real-Time Market}
\acrodef{stu}[STU]{Solar Thermal Unit}
\acrodef{csp}[CSP]{Concentrated Solar Plant}
\acrodef{ess}[ESS]{Energy Storage System}
\acrodef{bess}[BESS]{Battery Energy Storage System}
\acrodef{res}[RES]{Renewable Energy Source}
\acrodef{soc}[SoC]{State of Charge}
\acrodef{ndres}[ND-RES]{Non-dispatchable Renewable Energy Source}
\acrodef{dres}[D-RES]{Dispatchable Renewable Energy Source}
\acrodef{ro}[RO]{Robust Optimization}
\acrodef{sth}[S-TH]{Solar Thermal Unit}
\acrodef{vpp}[VPP]{Virtual Power Plant}
\acrodef{pv}[PV]{photovoltaic}
\acrodef{bes}[BES]{Battery Energy Storage}
\acrodef{phes}[PHES]{Pumped Hydro Energy Storage}
\acrodef{opf}[OPF]{Optimal Power Flow}
\acrodef{pcc}[PCC]{Points of Common Coupling}
\acrodef{dsm}[DSM]{Demand Side Management}
\acrodef{der}[DER]{Distributed Energy resource}

\newcommand{\add}[1]{\textcolor{black}{#1}}

\makeatletter
\let\old@ps@headings\ps@headings
\let\old@ps@IEEEtitlepagestyle\ps@IEEEtitlepagestyle
\def\psccfooter#1{%
    \def\ps@headings{%
        \old@ps@headings%
        \def\@oddfoot{\strut\hfill#1\hfill\strut}%
        \def\@evenfoot{\strut\hfill#1\hfill\strut}%
    }%
    \def\ps@IEEEtitlepagestyle{%
        \old@ps@IEEEtitlepagestyle%
        \def\@oddfoot{\strut\hfill#1\hfill\strut}%
        \def\@evenfoot{\strut\hfill#1\hfill\strut}%
    }%
    \ps@headings%
}
\makeatother

\begin{document}

\setlength{\abovedisplayskip}{0.5pt}
\setlength{\belowdisplayskip}{0.5pt}

\title{Modeling Demand Flexibility of RES-based Virtual Power Plants 
\vspace*{-0.25cm}}

\author{
\IEEEauthorblockN{Oluwaseun Oladimeji, {\'{A}}lvaro Ortega, Lukas Sigrist, Pedro S\'{a}nchez-Mart\'{i}n, Enrique Lobato, Luis Rouco}
\IEEEauthorblockA{Institute for Research in Technology, Comillas Pontifical University, Madrid, Spain}
\vspace*{-6mm}
}

\maketitle

\begin{abstract}
In this paper, an approach to evaluate the benefits of demand flexibility for \acp{vpp} is presented. The flexible demands chosen in this study are part of a renewable energy source-based \ac{vpp} that participates in \ac{dam} and \ac{idm} and has dispatchable and non-dispatchable assets. A demand model with bi-level flexibility is proposed: the first level is associated with \ac{dam}, whereas the second level is related to \ac{idm} sessions. Simulations are carried out considering a 12-node network to ascertain the eventual impacts of modeling demand flexibility on \ac{vpp} operation. \add{The market structure considered in the case study resembles the different trading floors in the Spanish electricity market. Results obtained show that the proposed demand flexibility scheme increases the overall profit of the \ac{vpp}, as well as the revenues of the demand owners without disrupting the consumer's comfort.}
\end{abstract}
\color{black}
\vspace*{-5mm}

\thanksto{\noindent This project has received funding from the European Union’s Horizon 2020 research and innovation programme under grant agreement No 883985 (POSYTYF project).}

\section{Introduction} \label{sec: intro}
\ac{dsm} for mutual benefits of consumer and the utility has been around for several decades. The definition, constituent parts and strategies for managing consumer's load have also evolved through this time \cite{gellings1985concept}. These strategies have included peak clippings, valley filling, load shifting, strategic conservation, strategic load growth, and flexible load shape, among others \cite{strbac2008demand}, \cite{palensky2011demand}. Today, popular loads that have been included in \ac{dsm} studies comprise electric vehicles (EVs) and large thermal-storage air-conditioning systems. These have been chosen largely due to their possession of storage and/or inertia provision capabilities. In many studies of load management however, a main challenge with \ac{dsm} models is the absence of appropriate measures of benefits for consumers that provide such flexibility actions.

In \cite{de2018local}, two types of load, flexible and non-flexible, were presented but assigning prices to this flexibility provision was not addressed. In a study with multiple \acp{der}, EVs serving as controllable electrical loads were used as a source of flexibility for load shifting. The objective here was increasing profitability of the \ac{der} assets \cite{Georgiev2016}. Similar to the previous study, authors of \cite{Hungerford2019} proposed the utilization of flexible loads to absorb the variable renewable generation during the day instead of simply doing load shifting from peak demand periods to times of low prices in the middle of the night. Moreover,  communication systems are increasingly permeating the power system with the usage of more Internet of Things devices. Thus, an information-rich energy system with flexible and responsive electrical loads that have storage capabilities can be built to respond to variable \acp{res} \cite{taneja2013flexible}.

\color{black}While these advancements are noteworthy, studies have given less attention to load owners having full control over what value they assign to this flexibility provision. 
This question can be coupled with seeking a solution for the market participation of stochastic nondispatchable \acp{res}. 
Aggregating \acp{res} within a \ac{vpp} is a solution to provide a better controlled output and to make \ac{res} more competitive in electricity markets \cite{marinescu2021dynamic}. 
In contrast to our present work, \acp{vpp} with flexible demands and who offer prices to those demand owners for their flexibility provision have not been taken into account in previous literature.

To fill this gap, a \ac{res}-based \ac{vpp} incorporating flexible demand for flexibility provision while participating in energy markets is presented.
\color{black}We propose a demand model with bi-level flexibility associated with different energy market trading sessions. In contrast to other studies, the demand owner maintains high level of control over its own consumption by setting different profiles which the \ac{vpp} manager can choose from in \ac{dam} at a cost. Moreover, the demand owner allows tolerances around the chosen profile at \ac{idm} at no cost. \color{black}
A network-constrained unit commitment model is formulated to submit \ac{vpp} \ac{dam} auctions and then subsequently \ac{idm} bids to correct for deviations from forecasts and other variations. \color{black}

In this study, we consider that the \ac{vpp} participates in the Spanish energy market. The energy market in Spain is well tailored for the studies carried out in this work and the authors have extensive knowledge of its market structure. The methodology proposed can be readily applied in other energy markets with similar structure and trading sessions.
\color{black}

\section{VPP Modelling} \label{sec: vpp model}
This section presents the formulation of the \ac{vpp} model used in this paper. The \ac{vpp} comprises flexible demands (industrial, airport and residential loads), \acf{dres} (hydro), \acfp{ndres} (wind power plant, solar PV) and \acf{stu} with storage capability. These assets are distributed across the electrical network and connected to the main grid through one or more \ac{pcc}. 
\color{black}The \ac{vpp} components and electrical network were presented in \cite{marinescu2021dynamic}.
\color{black}The business model considered for the \ac{vpp} is the maximization of its aggregated profit by optimally scheduling the generation and demand assets in its portfolio.

\add{The formulation for each asset class is discussed in the remainder of this section. \acp{dres} are modeled like conventional power plants \cite{baringo2020virtual} with linearized operation costs of the dispatchable assets. Network constraints are formulated by using a DC power flow approximation \cite{kargarian2016toward}. \acp{dres} and DC power flow constraints are well-known and not discussed here due to space limitations.} \color{black}In addition, the \ac{stu} model is adopted from \cite{oladimeji2022optimal}.
\add{The objective function of the \ac{vpp}, and constraints for the energy balance at the PCCs and the \acp{ndres} are presented in the following subsections.}
\color{black}

\subsection{Nomenclature}
This subsection presents the notation and nomenclature used in this section and in Section~\ref{sec: flex demands}. \vspace*{-8mm}
\begin{small}
\setlength{\nomitemsep}{0.05cm}

\nomenclature[I, 01]{$b \in \mathscr{B}/\mathscr{B}^m$}{Network Buses / Network buses Connected to Main Grid}
\nomenclature[I, 01]{$c \in \mathscr{C}/\mathscr{C}_b$}{\ac{dres} / \ac{dres} connected to bus $b$}
\nomenclature[I, 01]{$d \in \mathscr{D}/\mathscr{D}_b$}{Demand / Demand connected to bus $b$}
\nomenclature[I, 02]{$\ell \in \mathscr{L}$}{Network lines}
\nomenclature[I, 01]{$k \in \mathscr{K}$}{\ac{idm} sessions}
\nomenclature[I, 03]{$p \in \mathscr{P}$}{Demand profiles}
\nomenclature[I, 03]{$r \in \mathscr{R}/\mathscr{R}_b$}{\ac{ndres} / \ac{ndres} connected to bus $b$}
\nomenclature[I, 03]{$t \in \mathscr{T}$}{Time periods}
\nomenclature[I, 04]{$\theta \in \Theta/\Theta_b$}{\ac{stu} / \ac{stu} connected to bus $b$ \vspace{5pt}}

\nomenclature[P, 01]{$C_c^0/C_c^1$}{Shut-down/start-up cost of \acp{dres}\nomunit{[\euro]}}  
\nomenclature[P, 01]{$C_c^{\rm V}$}{Variable production cost of \acp{dres}\nomunit{[\euro/MWh]}}

\nomenclature[P, 01]{$C_{d,p}$}{Cost of load profile $p$ of demand\nomunit{[\euro]}}
\nomenclature[P, 02]{$\ubar{E}_d$}{Min energy consumption of demand $d$ throughout the planning horizon\nomunit{[MWh]}}
\nomenclature[P, 03]{$P_{d,p,t}$}{Max hourly consumption of profile $p$ of demand $d$ \nomunit{[MW]}}
\nomenclature[P, 03]{$\ubar{P}_{d,t}/\bar{P}_{d,t}$}{Lower/upper bound variations of the power consumption of demand $d$ in time $t$\nomunit{[\%]}}
\nomenclature[P, 05]{$\ubar{R}_d/\bar{R}_d$}{Down/up ramping limit of demand $d$\nomunit{[MW/h]}}

\nomenclature[P, 03]{$\ubar{P}_{r,t}$}{Min production of \ac{ndres} in time $t$\nomunit{[MW]}}
\nomenclature[P, 04]{$\check{P}_{r,t}$}{Available production of \ac{ndres} in time $t$\nomunit{[MW]}}


\nomenclature[P, 05]{$\bar{P}_b^m$}{Maximum power that can be traded with the main grid at bus $b$ \nomunit{[MW]}}

\nomenclature[P, 07]{$\lambda^{\rm DA}_t$}{\ac{dam} price in time $t$\nomunit{[\euro/MWh]}}
\nomenclature[P, 07]{$\lambda^{\rm ID}_{k,t}$}{Price of \ac{idm} session $k$ in time  $t$ \nomunit{[\euro/MWh]} \vspace{5pt}}
\nomenclature[P, 07]{$\Delta t$}{Duration of time periods\nomunit{[h, min]}}

\nomenclature[V, 01]{$p_{c,t}$}{Power generation of \acp{dres} in time $t$\nomunit{[MW]}}

\nomenclature[V, 01]{$p_{d,t}$}{Power consumption of demand in time $t$\nomunit{[MW]}}
\nomenclature[V, 01]{$u_{d,p}$}{Binary variable to select demand profile\nomunit{[$0/1$]}}

\nomenclature[V, 01]{$p_{r,t}$}{Power generation of \ac{ndres} in time $t$\nomunit{[MW]}}

\nomenclature[V, 01]{$p_{\theta,t}$}{Electrical power generation of \ac{stu} in time $t$\nomunit{[MW]}}

\nomenclature[V, 01]{$p_{\ell,t}$}{Power flow through network of line $\ell$ in time $t$\nomunit{[MW]}}

\nomenclature[V, 01]{$p_{b,t}^m$}{Power scheduled to be bought from/sold to the \ac{dam} and \ac{idm} markets at bus $b$ in time  $t$\nomunit{[MW]}}

\nomenclature[V, 01]{$p_t^{\rm DA}/p_{k,t}^{\rm ID}$}{Total power traded in the \ac{dam}/\ac{idm} in time $t$\nomunit{[MW]}}

\renewcommand{\nomname}{}
\printnomenclature[1.5cm]
\end{small} 

\subsection{Profit Maximization Objective}
Due to small volumes of energy traded in the \ac{idm} relative to \ac{dam} and modest price differences between these markets in the Spanish system \cite{chaves2015spanish}, the objective functions in \ac{dam} and \ac{idm} are decoupled in this work.
Each \ac{idm} further has  associated constraints to cater for updates or changes in forecasts of stochastic sources.
In \ac{dam}, the objective function \eqref{obj:DAM_det} is the maximization of the obtainable profit by the \ac{vpp} assets calculated as the revenue from power trades minus cost of operating \ac{dres} and cost of selecting a particular load profile. Operation cost of \ac{ndres} are not considered in the objective function due to their relative low value when compared to costs of dispatchable sources.
For the different \ac{idm} sessions, the benefit \eqref{obj:IDM_det} is calculated over changes in traded power $\Delta p_{k,c,t}$ between: (i) \ac{dam} and first \ac{idm} trading period and (ii) other subsequent \ac{idm} sessions. Cost of choosing a specific load profile is not included while computing objective of \ac{idm} because this choice is previously made during \ac{dam} participation and must be accounted for only once.
\begin{multline} \label{obj:DAM_det}
    \max_{\Xi^{\rm DAM}} 
    \sum_{t\in\mathscr{T}} 
    \left[ 
    \lambda^{\rm DA}_t p_t^{\rm DA} \Delta t - 
    \sum_{c\in\mathscr{C}} 
    \left( 
    C_c^{\rm V} p_{c,t} \Delta t + 
    c_{c,t}^0 + 
    c_{c,t}^1 
    \right) 
    \right] \\ 
    - \sum_{d \in \mathscr{D}} 
    \sum_{p \in \mathscr{P}} 
    C_{d,p}u_{d,p} \qquad \qquad \qquad \qquad \qquad \qquad \quad \vspace{-1mm}
\end{multline}
\begin{multline} \label{obj:IDM_det}
    \max_{\Xi^{\rm IDM}_k} 
    \sum_{t=\tau_k}^{\left| \mathscr{T} \right|} 
    \Bigg[
    \lambda^{\rm ID}_{k,t}  
    p_{k,t}^{\rm ID} \Delta t \\
    - \sum_{c\in\mathscr{C}} 
    \left( 
    C_c^{\rm V} \Delta p_{k,c,t} \Delta t + 
    c_{k,c,t}^0 + 
    c_{k,c,t}^1 
    \right) 
    \Bigg], 
    \quad 
    \forall k \in \mathscr{K}
\end{multline}

\subsection{Energy Balance}
Energy balance constraints common to both market stages are modeled in \eqref{cons:balance common} whereas those specific to \ac{dam} and \ac{idm} are formulated in \eqref{cons:balance_pDA2} and \eqref{cons: energy_bal_idm3} respectively.
Equation \eqref{cons:balance_mg} gives energy balance at the \ac{pcc} with the main grid, while \eqref{cons:balance_no_mg} is the balance for all other buses in the \ac{vpp} network at every time period \cite{van2014dc}. The difference between both equations is presence of $p_{b,t}^m$ at the main grid representing scheduled power to be traded with other market participants. This power available for trading (buy or sell) is set within prespecified bounds in \eqref{cons: trade_bound}. \vspace*{-2mm}
\begin{subequations} \label{cons:balance common}
\begin{multline} \label{cons:balance_mg}
    \sum_{c\in\mathscr{C}_b} p_{c,t} \!+\!\! 
    \sum_{r\in\mathscr{R}_b} p_{r,t} \!+\!\! 
    \sum_{\theta\in\Theta_b} p_{\theta,t} \!-\!\! \!\!\!
    \sum_{\ell | i(\ell)=b} \!\!p_{\ell,t} +\! \!\!\!
    \sum_{\ell | j(\ell)=b} \!\!\!p_{\ell,t} \\
    = p_{b,t}^m \!+\!\! 
    \sum_{d\in\mathscr{D}_b} p_{d,t}~, 
    \qquad 
    \forall b \in \mathscr{B}^m, 
    \forall t \in \mathscr{T} 
\end{multline}
\begin{multline} \label{cons:balance_no_mg}
    \sum_{c\in\mathscr{C}_b} p_{c,t} \!\!+\!\! 
    \sum_{r\in\mathscr{R}_b} p_{r,t} \!\!+\!\! 
    \sum_{\theta\in\Theta_b} p_{\theta,t} -\!\! \!\!
    \sum_{\ell | i(\ell)=b} \!\!p_{\ell,t} +\!\!  \!\!
    \sum_{\ell | j(\ell)=b} \!\!p_{\ell,t} \\ 
    = \sum_{d\in\mathscr{D}_b} p_{d,t}~,
    \qquad 
    \forall b \in \mathscr{B} \setminus \mathscr{B}^m, 
    \forall t\in \mathscr{T}
\end{multline} \vspace{1ex}
\begin{IEEEeqnarray}{lr}
    -\bar{P}_b^m 
    \leq p_{b,t}^m 
    \leq \bar{P}_b^m ~,
    & \qquad \qquad \;\;
    \forall b \in \mathscr{B}^m, 
    \forall t \in \mathscr{T}  \label{cons: trade_bound}
\end{IEEEeqnarray} \vspace{0.0em}
\end{subequations} 
\vspace{-2mm}

\subsubsection{DAM Formulation} 
Equation \eqref{cons:balance_pDA2} ensures that summation of traded power at all \acp{pcc} is equivalent to the total power available for trading by \ac{vpp}.
\begin{IEEEeqnarray}{lr} 
    p_t^{\rm DA} = \!\!
    \sum_{b\in \mathscr{B}^m} p_{b,t}^m ~, 
    & \qquad \qquad \qquad \qquad \qquad
    \forall t \in \mathscr{T} \label{cons:balance_pDA2}
\end{IEEEeqnarray}

\subsubsection{IDM Formulation}
For \ac{idm} sessions, the \ac{idm} offers/bids do not substitute those submitted in the \ac{dam}, but rather, they are \textit{adjustments} of the \ac{dam} offers/bids as reflected in \eqref{cons: energy_bal_idm3}.
\color{black}The rationale behind such adjustments can be due to unplanned maintenance of generators, changes in \ac{ndres} outputs, demand changes and/or line faults.
\color{black}
\begin{IEEEeqnarray}{lr}
    p_t^{\rm DA^*} \!\!+\!  
    \sum_{\kappa=1}^{k-1}p_{\kappa,t}^{\rm ID^*} \!+\!  
    p_{k,t}^{\rm ID} = \!\!\!
    \sum_{b\in \mathscr{B}^m} p_{b,t}^m ~,
    & \qquad
    \forall k \in \mathscr{K}, 
    \forall t \geq \tau \qquad \label{cons: energy_bal_idm3}  \vspace{1ex}
\end{IEEEeqnarray} 
where $p_t^{\rm DA^*}\!\!$ and $p_t^{\rm ID^*}\!\!$ are solutions of the \ac{dam} and previous \acp{idm} respectively.
Note that the only difference with respect to \eqref{cons:balance common} is
the time index, $\forall t \in \mathscr{T}$, which is replaced with $\forall t \geq \tau$, where $\tau$ is the first delivery period of current \ac{idm} session.

\subsection{\aclp{ndres}} \label{sec: ndres}
The \acp{ndres} modeled in \eqref{cons:ndres} comprise mainly wind power and solar PV plants. The lower bound represents the asset technical minimum (e.g., cut-in speed for wind plant) while the output is bounded above by the available stochastic source.
\vspace*{-3mm}
\begin{IEEEeqnarray}{lr} \label{cons:ndres}
    \ubar{P}_{r,t}
    \leq 
    p_{r,t}
    \leq 
    \check{P}_{r,t}~,
    & \qquad \qquad \qquad \qquad
    \forall r \in \mathscr{R}, 
    \forall t \in \mathscr{T} \qquad
\end{IEEEeqnarray}

\section{Flexible Demands} \label{sec: flex demands}
This section presents the model of the demand flexibility in a \ac{vpp} participating in energy markets that is proposed in this work.
The model comprises two levels of flexibility, each associated with \ac{dam} and \ac{idm} market sessions respectively. \color{black}

\subsection{DAM Formulation} \label{ssec: flex DAM}
An optimal load profile is selected during \ac{dam} as the first flexibility level. For each demand $d$, \eqref{cons: demand profiles} and \eqref{cons: demand choice} guarantees the choice of a single profile $p$ out of all the available ones. Prior to market participation, various profiles are prepared by the demand owners/aggregators and communicated to the \ac{vpp} manager. A particular profile might serve as the \textit{default}, i.e., the consumption profile that will be followed by the demand if no market participation is considered. Additional profiles which the load owner can follow but might be operationally costly can be presented. However, compensation will be required from the \ac{vpp} if those are chosen. Take an instance of residential demands; the default profile might feature double peaks at 09:00 and 20:00 whereas another profile features a shift of these load peaks to 07:00 and 21:00 respectively. If the second profile is selected as optimal by the \ac{vpp}, the cost of \textit{operational inconvenience} must be paid to the load owner.
\vspace*{-3mm}
\begin{subequations} \label{cons:flexload_DAM}
\begin{IEEEeqnarray}{lr}
    p_{d,t} = 
    \sum_{p \in \mathscr{P}} 
    P_{d,p,t} u_{d,p} ~,
    & \qquad \qquad \qquad
    \forall d \in \mathscr{D}, 
    \forall t \in \mathscr{T} \qquad \; \label{cons: demand profiles} \\
    \sum_{p \in \mathscr{P}} u_{d,p} = 
    1 ~,
    & 
    \qquad
    \forall d \in \mathscr{D} 
    \phantom{,\forall p \in \mathscr{P}} \qquad \label{cons: demand choice} 
\end{IEEEeqnarray} \vspace{-2mm}
\end{subequations}

\subsection{IDM Formulation}
The second level of demand flexibility is provided during the different \ac{idm} sessions, formulated in \eqref{cons:flexload_IDM}.  At \ac{idm}, the load profile selected from \ac{dam} cannot be changed. However, \color{black}the demand owner allows the VPP manager to vary the consumption within predefined threshold (symmetric or not) around that selected profile ($P_{d,p,t}^*$) as presented in \eqref{cons: load flex}. \color{black} \color{black}Equations~\eqref{cons: load ramp1} and \eqref{cons: load ramp2} define the ramps of the demand profile from one period to the next. Finally, \eqref{cons: load minimum} ensures that, over the total duration of the current \ac{idm} session plus the periods covered in previous sessions, a minimum amount of energy is consumed. The energy values settled in previous periods, $p_{d,t}^*$, are thus accounted for during subsequent \ac{idm}.
\begin{subequations} \label{cons:flexload_IDM}
\begin{multline} \label{cons: load flex}
    \left(
    1 - \ubar{P}_{d,t}
    \right) 
    P_{d,p,t}^{*}
    \leq p_{d,t} 
    \\
    \qquad \qquad \leq   
    \left(
    1 + \bar{P}_{d,t} 
    \right) 
    P_{d,p,t}^{*}~,\!
    \qquad 
    \forall d \in \mathscr{D}, 
    \forall t \geq \tau \,\, \hspace{-0.6em}
\end{multline}
\begin{IEEEeqnarray}{lr}
    p_{d,t} - 
    p_{d,(t-1)} 
    \leq 
    \bar{R}_d \Delta t~,
    & 
    \forall d \in \mathscr{D}, 
    \forall t \geq \tau \quad \label{cons: load ramp1} \\
    p_{d,(t-1)} - 
    p_{d,t} 
    \leq 
    \ubar{R}_d \Delta t~,
    & 
    \forall d \in \mathscr{D}, 
    \forall t \geq \tau \quad \label{cons: load ramp2} \\
    \ubar{E}_d 
    \leq 
    \sum_{t=1}^{\tau-1} p_{d,t}^* \Delta t + 
    \sum_{t= \tau}^{\left| \mathscr{T} \right|} p_{d,t} \Delta t~,
    & \qquad
    \forall d \in \mathscr{D} 
    \phantom{,\forall t \geq \tau} \quad \label{cons: load minimum}
\end{IEEEeqnarray}
\end{subequations}

\color{black}

\section{Case Study} \label{sec: case study}
This section presents the case studies considered to test and validate the RES-based \ac{vpp} model proposed in this paper. Section \ref{subsec: vpp description} outlines the \ac{vpp} topology considered, which resembles a subarea of the southern region of the Spanish grid. The input data that for the model is described in Section~\ref{subsec: input params}.

\subsection{VPP Description} \label{subsec: vpp description}
\begin{figure}[!htbp]
    \centering
    \vspace*{-3mm}
    \includegraphics[width=\linewidth, height=1.6in]{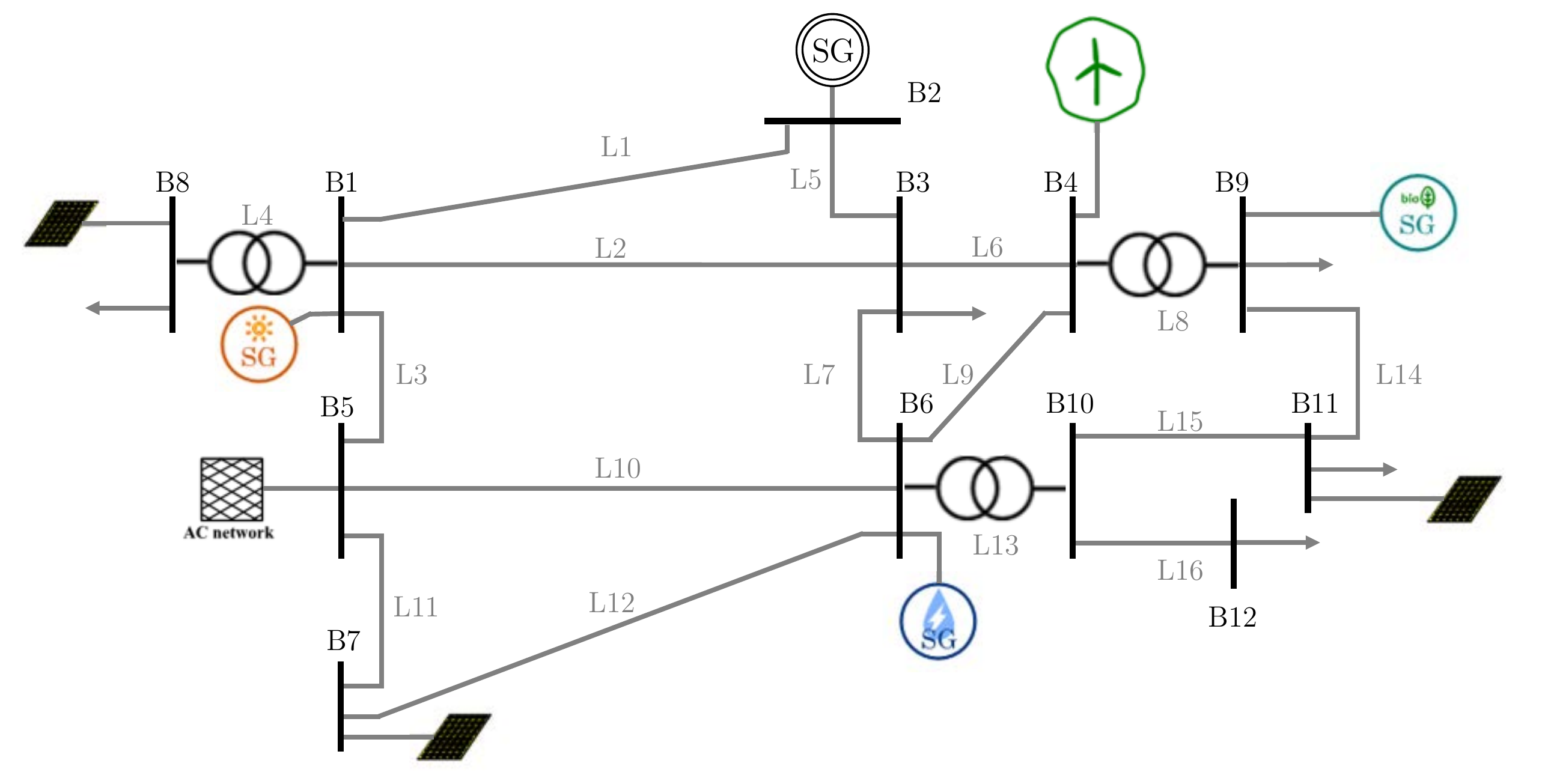}
    \caption{12-node network for test cases}
    \vspace*{-4mm}
    \label{fig:network structure}
\end{figure}

The \ac{vpp} assets are distributed across a 12-node network connected to a main grid through the \ac{pcc} (bus 5) as shown in Fig.~\ref{fig:network structure}. 
The demands considered are industrial, airport and residential loads (buses 3, 9 and 12)  with minimum daily consumption of $800$, $580$, and $600$ MWh respectively. Three profiles are associated with each demand and total consumption for each profile is the same.
Capacity of hydro (bus 6) is $111$ MW. Wind power plant (bus 4), solar PV (bus 8) and \ac{stu} (bus 1) each have capacities of $50$ MW. 

\subsection{Input Parameters} \label{subsec: input params}
\begin{figure*}[!htbp]
    \centering
    \includegraphics[width=\linewidth]{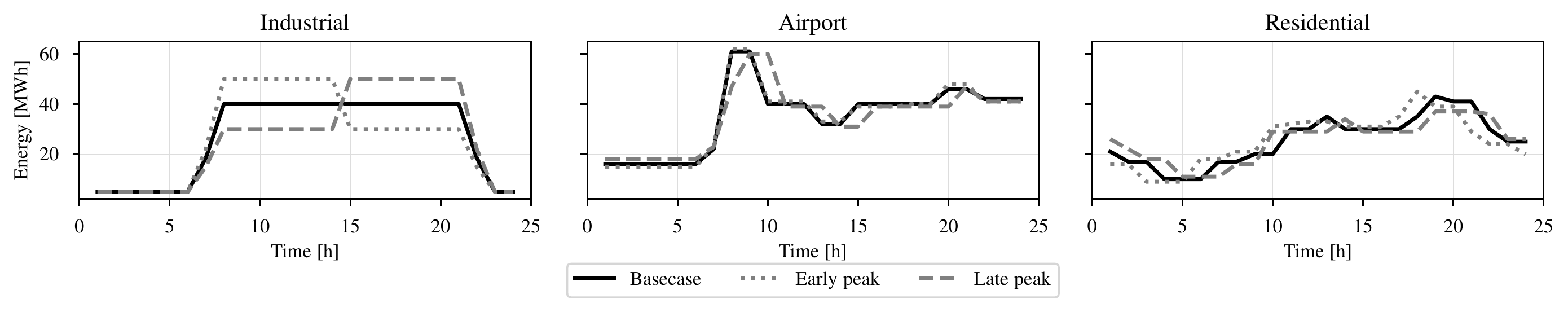}
    \caption{Demand profiles}
    \label{fig:load profiles}
\end{figure*}

A time horizon of $24$ hours with an hourly resolution is utilized during \ac{dam}, while a subset of the $24$ hours is used for each \ac{idm} session (see~\cite{chaves2015spanish}). 
Figure.~\ref{fig:load profiles} shows the three profiles considered for each of the loads. The basecase profile simulates \textit{default activities} of electricity consumers (see Section~\ref{ssec: flex DAM}).
The other two profiles (early and late peaks) are designed to perform load shifting around the basecase. Other profiles can also be easily added to capture other modes of demand interaction with the market. During \ac{idm}, the demand owner allows a percentage of tolerance for demand movement over the selected profile at \ac{dam}. In this case study, symmetric tolerances from $10-50\%$ over the total demand are considered. Note that the highest tolerance values considered serve only as an illustrative study, as likely no demand owner would allow such variations in a short notice (sometimes as short as 1 hour).

\vspace*{-3mm}
\section{Results} \label{sec: results}
To test the effectiveness of our model using renewable sources, two distinct operation days are identified: a clear sunny day and a day with intermittent cloud covers. 
Two aspects are thoroughly discussed for each day, namely: a) optimal price to be offered by the \ac{vpp} to the demand owners and \color{black} b) effect of the flexibility provided by demand on the profits accrued and operation of the \ac{vpp}.

\subsection{Optimal Price Offered by VPP to Demand Owners}
\color{black}First, we provide a solution to the benefits for demand owners and show that there is a maximum daily total ($C_{d,p}$) and a price per MWh that the \ac{vpp} is prepared to pay demand owners such that it is still profitable. \color{black}The default load profile has zero cost because the demand owner follows it regardless of other events. \ac{vpp} then decides what price to pay demand owners for other profiles such that it's benefits are not eroded. 

\subsubsection{Clear day}
The optimal choices of demand profiles obtained for \ac{vpp} operation on a clear day is shown in Table~\ref{table:demand profile choice clear day}. \color{black} These are obtained when $C_{d,p}$ for every other profile is set to zero. 
\color{black} Industrial and residential loads have the early peak and late peak as optimal profiles respectively. 
\color{black} To determine the cut-off costs before the basecase is selected, costs for other demand profiles are then gradually increased (simultaneously for all demands). It was then observed that \color{black}
the \ac{vpp} manager is only willing to pay up to \euro$320$/day to the industrial load owner.
The availability of the default profile ($C_{d,p}=0$) acts as an economic viability check for the \ac{vpp}. The \euro$320$ is thus the maximum price offering for the early peak industrial load profile, after which it is no more profitable for \ac{vpp}. The industrial late peak profile is not profitable at all and neither is the residential early peak profile.
In the case of the airport load, there are more options for both the airport demand owner and the \ac{vpp} manager. The early peak profile is profitable for the \ac{vpp} manager as a cost up until \euro$500$/day while the late peak profile is only profitable until \euro$305$/day. 

As observed in Fig.~\ref{fig:load profiles}, the different airport profiles have only subtle differences. Additionally, the load profiles have a benefit-inducing relationship with the evolution of the market price on this clear day. Indeed, the first profile peaks (periods 8-10) correspond to lower electricity prices and thus account for the higher costs which the \ac{vpp} is prepared to pay.
\color{black}When the demand cost per MWh of the optimal profiles are computed, proof of the hypothesis is further demonstrated. For the airport profiles, $C_{d,p}$ = \euro$10.42/$MWh for $48$ MWh shifted and \euro$7.44/$MWh for $41$ MWh shifted of the early and late peak profiles respectively. Contrast this with industrial early peak profile where $C_{d,p}$ = \euro$4.35/$MWh for $73.5$ MWh shifted. 
Thus, the \euro/day offer from VPP depends a little less on the magnitude of energy shifted and more on the load profile shape with respect to market price. 
Finally, these costs per MWh associated with \ac{dam} represent lower bounds on the price that \ac{vpp} is willing to pay demand owners. With the tolerances allowed in \ac{idm}, \ac{vpp} benefits increase and the prices paid to demand owners are expected to follow the same trend. 
\color{black}
\begin{table}[!tbp]
  \centering
  \caption{Demand profile choice at non-zero costs on a clear day}
 \vspace{1mm}
  \begin{tabular}{llcc}
    \toprule
    \multicolumn{1}{c}{\textbf{Demand}} & \multicolumn{1}{c}{\textbf{Basecase}} &    \multicolumn{1}{c}{\textbf{Early peak}} & \multicolumn{1}{c}{\textbf{Late peak}} \\
    \midrule
    \multirow{1}{*}{Industrial} & chosen when & \multirow{1}{*}{optimal} & \multirow{1}{*}{\xmark} \\
    & cost$>$€320/day &  & \\[0.25em]
    \multirow{1}{*}{Residential} & chosen when & \multirow{1}{*}{\xmark} & \multirow{1}{*}{optimal} \\
    & cost$>$€180/day &  & \\[0.25em]
    \multirow{1}{*}{Airport} & chosen when   & \multirow{1}{*}{optimal} & \multicolumn{1}{l}{suboptimal - chosen}  \\
    & cost$>$€500/day  & & \multicolumn{1}{l}{when cost$>$€305/day}  \\    
    \bottomrule
  \end{tabular}
  \label{table:demand profile choice clear day}
\end{table}

\subsubsection{Cloudy day}
On a cloudy day, it is observed that for all demands, the late peak profile is not optimal for the \ac{vpp} operation. Additionally, the \ac{vpp} manager is unwilling to pay amounts as high as the clear day case. Only \euro$90$/day and \euro$260$/day is offered to the industrial and airport load respectively. For the residential load, the deviations from the default load are not optimal at all and the \ac{vpp} manager thus, offers no price to the residential load owner.  \vspace*{-4mm}
\begin{table}[!htbp]
  \centering
  \caption{Demand profile choice at non-zero costs on a cloudy day}
 \vspace{1mm}
  \begin{tabular}{llcc}
    \toprule
    \multicolumn{1}{c}{\textbf{Demand}} & \multicolumn{1}{c}{\textbf{Basecase}} &    \multicolumn{1}{c}{\textbf{Early peak}} & \multicolumn{1}{c}{\textbf{Late peak}} \\
    \midrule
    \multirow{1}{*}{Industrial} & chosen when & \multirow{1}{*}{optimal} & \multirow{1}{*}{\xmark} \\
    & cost$>$€90/day &  & \\[0.25em]
    \multirow{1}{*}{Residential} & optimal & \multirow{1}{*}{\xmark} & \multirow{1}{*}{\xmark} \\
    & &  & \\[0.15em]
    \multirow{1}{*}{Airport} & chosen when   & \multirow{1}{*}{optimal} & \multirow{1}{*}{\phantom{suboptim} \xmark \phantom{chosen wh}}  \\
    & cost$>$€260/day  &  &  \\    
    \bottomrule
  \end{tabular}
  \label{table:demand profile choice cloudy day}
\end{table}  \vspace*{-5mm}

\color{black}
\subsection{Effects of Demand Flexibility on VPP Profits}
This section discusses: 1) effects of demand flexibility on profitability of the \ac{vpp} and 2) effects of percentage flexibility allowed by the demand owners.

\begin{figure}[!htbp]
     \centering
     \vspace*{-4mm}
     \includegraphics[width=1\linewidth]{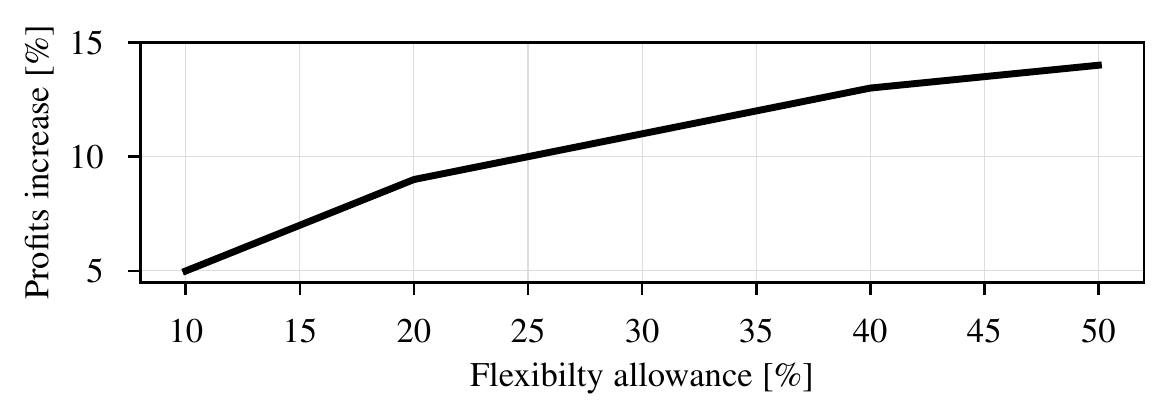}
     \caption{Demand profile flexibility on \ac{vpp} profits on a clear day}
     \vspace*{-1mm}
     \label{fig:Dem prof flex on vpp profits clear day}
\end{figure}

On a clear day, the \ac{vpp} manager can realize a profit of \euro$68259$ from participating in the \ac{dam} (first demand flexibility level). However, at the end of the operation horizon (i.e. final \ac{idm} session), this profit can increase by up to $14\%$ attributable to the second demand flexibility level alone beyond its \ac{dam} objective. The relationship between second-level flexibility allowance and the associated effect on profits is shown in Fig.~\ref{fig:Dem prof flex on vpp profits clear day}. \color{black}Beyond 40\%, profits do not show any considerable additional increase and any allocation on load movement beyond this percentage only leads to more saturation on the profits curve. 

\begin{figure}[!htbp]
     \centering
     \vspace*{-3mm}
     \includegraphics[width=1\linewidth]{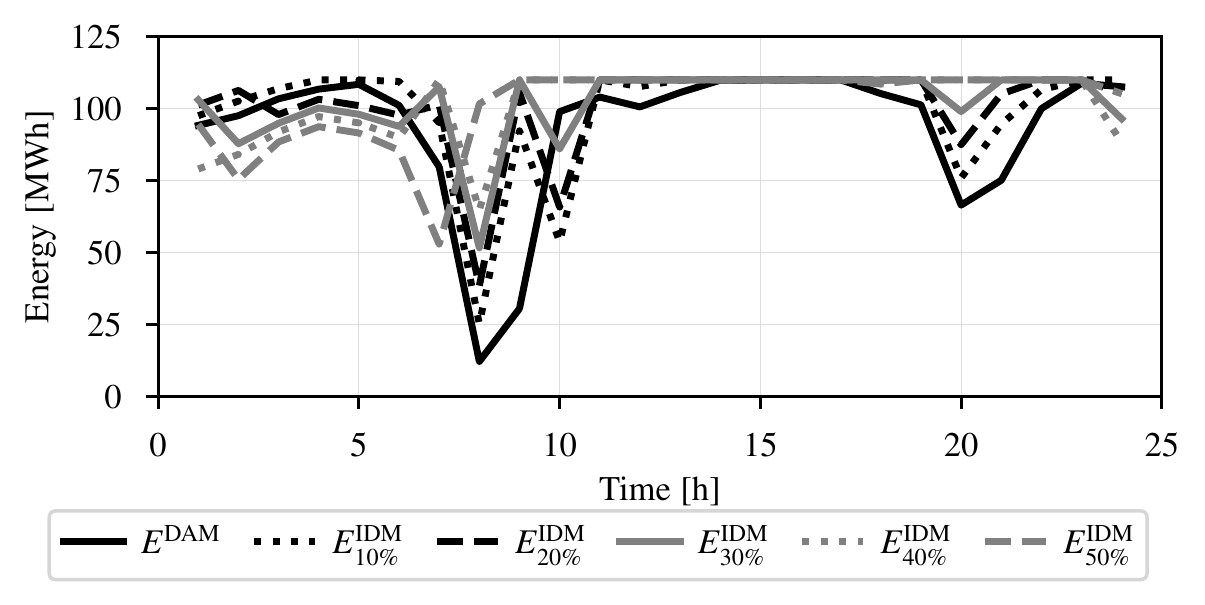}
     \caption{Scheduled and final traded power based on demand flexibility}
     \vspace*{-1mm}
     \label{fig:traded power clear day}
\end{figure}

Figure.~\ref{fig:traded power clear day} shows the impact of demand flexibility on total power traded in \ac{dam} and total traded after all \ac{idm} sessions. \color{black}It is observed that increasing flexibility provision reduces the valley in \ac{dam} offer around the instant h = 8 and leads to an increase in offers in the energy market while at the same time flattening the curve in the straight section from periods 9 to 19 but keeping total energy consumption constant. \color{black}

\begin{figure*}[!b]
     \centering
     \vspace*{0mm}
     \includegraphics[width=1\linewidth]{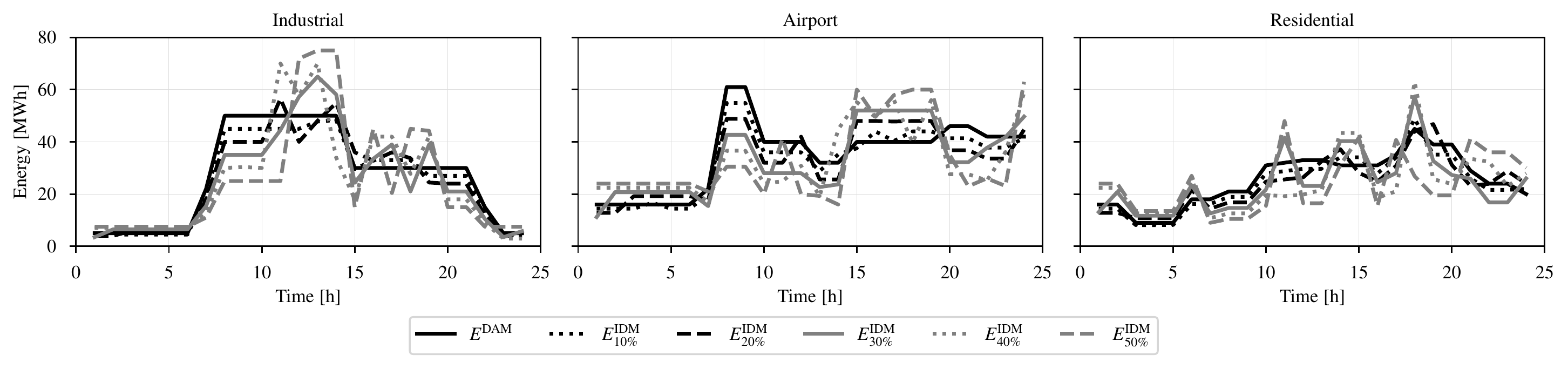}
     \caption{Scheduled and final demand outputs on a clear day}
     \label{fig:final dem output on clear Day}
\end{figure*}

Final demand output relative to \ac{dam} offer is shown in Fig.~\ref{fig:final dem output on clear Day}. Up to 20\% flexibility allowance on demands, the variations are smooth and would not necessarily alter normal system operation. However, from 30\% and upwards flexibility, higher consumption pattern is observed in demand as evidenced in the industrial and residential demands at periods 10-13 and 17-19. Increasing flexibility allowance gives rise to more of this peaks that might affect system performance. Demand owners would likely not allow such percentages anyway and we simply show a proof of concept here. \color{black}

\section{Conclusion} \label{sec: conclusion}
This paper presents a model for evaluating the impacts of flexibility provision by demands that are part of an \ac{res}-based \ac{vpp} that participates in energy markets.
The flexible demand model proposed has two levels of flexibility associated with different market sessions. The business model of the \ac{vpp} is such that it maximises its profit by dispatching its generators, carrying out self supply of the demands within its portfolio and utilizing the flexibility actions provided by said demands. 

Case studies were then analysed to determine the optimal price the \ac{vpp} is ready to offer demand owners. Impact of the flexibility provided on operation of \ac{vpp} and profits accrued at the end of the operation day is also discussed. Based on the studies carried out, it is concluded that there is a maximum price offering from \ac{vpp} to demand owners for each demand profile. Beyond this price, the default load will be selected. A sensitivity analysis showed that beyond 40\% flexibility allowance for the demands studied, there is a saturation of the profit and added flexibility might not be anymore profitable.

\bibliographystyle{IEEEtran}
\bibliography{refs}

\end{document}